\documentclass[12pt]{article}
\usepackage{amsfonts, amsmath, amssymb,amsthm,bbm,comment,fancyhdr,float,latexcad,makeidx,mathrsfs,verbatim}
\normalsize
\newcommand{\Ga}{\Gamma}

\newcommand{\La}{\Lambda}

\newcommand{\seq}{\subseteq}

\newcommand{\vs}{\vspace*}

\newcommand{\nin}{\noindent}

\newtheorem{mthm}{Theorem}[section]
\newtheorem{mylem}[mthm]{Lemma}
\newtheorem{myprn}[mthm]{Proposition}
\newtheorem{mycor}[mthm]{Corollary}
\newtheorem{mydef}[mthm]{Definition}
\newtheorem{myrem}[mthm]{Remark}
\newtheorem{mycon}[mthm]{Construction}
\newtheorem{myeg} [mthm]{Example}
\newtheorem{myque} [mthm]{Question}
\newenvironment{thm}{\begin{mthm}}{\end{mthm}}
\newenvironment{lem}{\begin{mylem}}{\end{mylem}}
\newenvironment{cor}{\begin{mycor}}{\end{mycor}}

\newenvironment{mdef}{\begin{mydef}\rm}{\end{mydef}}
\newenvironment{rem}{\begin{myrem}\rm}{\end{myrem}}

\newenvironment{ex}{\begin{myeg}\rm}{\end{myeg}}
\newenvironment{prof}{\noindent $Proof.$ \rm}{\hfill $\Box$}
\newenvironment{que}{\begin{myque}\rm}{\end{myque}}

\def \nin {\noindent}

\def \Lemma #1 {\vs{3mm}\nin {\bf Lemma #1} \it}
\def \Prop #1 {\vs{3mm}\nin {\bf Proposition #1} \it}
\def \Th #1 {\vs{3mm}\nin {\bf Theorem #1} \it}
\def \Cor #1 {\vs{3mm}\nin {\bf Corollary #1} \it}
\def \Ex #1 {\vs{3mm}\nin {\bf Example #1} \it}

\def \part #1 {\hfil\break\hglue 12pt {\rm (#1)~}}

\def\fs{\footnotesize}

\begin{document}
\title{\Large\bf On realizing zero-divisor graphs of po-semirings \thanks{The first author  is partly supported by Shanghai Jiaotong University Innovation Fund for Postgraduates $($NO. AE071202$)$.}}
\author{Houyi Yu\thanks{ yhy178@163.com (H.Y. Yu)}\,\, and \, Tongsuo Wu\thanks{Corresponding author,\, tswu@sjtu.edu.cn (T.S. Wu)}\\
{\small  Department of Mathematics, Shanghai Jiaotong University, Shanghai, $200240$, China}
}

\date{}
\maketitle

\begin{center}
\begin{minipage}{12cm}
\nin{\bf Abstract.} {\fs In this paper, we determine bipartite graphs and complete graphs with horns, which are realizable as zero-divisor graphs of po-semirings. As applications, we classify commutative rings $R$ whose annihilating-ideal graph $\mathbb {AG}(R)$ are either bipartite graphs or complete graphs with horns.}\\

\nin{\bf Key Words:} {\fs Po-semirings;  Graph properties; Zero-divisors; Annihilating-ideals of a commutative ring}\\

\nin{\bf MSC(2000):} {\fs 13A15; 05C75.}

\end{minipage}
\end{center}

\section{Introduction}
Throughout this paper, all rings are assumed to be commutative with identity. For a ring $R$,  let $Z(R)$ be its set of zero-divisors. The {\em zero-divisor graph} of $R$, denoted by $\Ga(R)$,  is a simple graph  (i.e., an undirected graph without loops and multiple edges) with vertices $Z(R)^*=Z(R)\backslash\{0\}$, such that distinct vertices $x$ and $y$ are adjacent if and only if $xy=0$. The concept for a ring was first introduced and studied by Beck in \cite{BECK} and further investigated by many authors, see  \cite{DAMN,AL,AMY}. Later, DeMeyer, McKenzie and Schneider \cite{{DMS}} extended the notion to commutative semigroups $S$ with $0$ in a similar manner. The idea establishes a connection between graph theory and algebraic theory and will be beneficial for those two branches of mathematics.

 For a ring $R$,  let $\mathbb I(R)$ be the set of ideals of $R$, $\mathbb A(R)$ the set of annihilating-ideals of $R$, where a nonzero ideal $I$ of $R$ is called an {\em annihilating-ideal} if there exists a non-zero ideal $J$ of $R$ such that $IJ=0$.  The {\em annihilating-ideal graph} $\mathbb {AG}(R)$ of $R$, first introduced and studied in \cite{BR}, is a graph with vertex set $\mathbb {A}(R)^*=\mathbb {A}(R)\backslash\{0\}$, such that distinct vertices $I$ and $J$ are adjacent if and only if $IJ=0$. The graph provides an excellent setting for studying some aspects of algebraic property of a commutative ring, especially, the ideal structure of a ring. Clearly, the graph $\mathbb {AG}(R)$ is an empty graph if and only if $R$ is an integral domain.

In fact, $\mathbb I(R)$ admits a natural algebraic structure, called a {\it po-semiring} by Wu, Lu and Li \cite{Wll}. Recall that a {\it commutative semiring} is a set $A$ which contains at least two elements $0,1$ and is equipped with two binary operations, $+$ and $\cdot$, called addition and multiplication respectively, such that the following conditions hold:

(1)  $(A, +,0)$ is a commutative monoid with zero element $0$.

(2)  $(A, \cdot,1)$ is a commutative monoid with identity element 1.

(3)  Multiplication distributes over addition.

(4)  $0$ annihilates $A$ with respect to multiplication, i.e., $0¡¤a = 0, \,\forall a \in A$.

Recall the following definition from \cite{Wll}.

\vs{3mm}\begin{mdef} \label{po-semiring}
A  {\em partially-ordered semiring} is a commutative semiring $(A, +, \cdot,0,1 )$, together with a compatible partial order $\le$,  i.e.,  a partial order $\le$ on the underlying set $A$ that is compatible with the semiring operations in the sense that it satisfies the following conditions:

(5) $x\le y$ implies $x+z\le y+z$, and

(6) $0\le x$ and $y\le z$  imply that $xy\le xz$
for all $x,y,z$ in $A$.

If $A$ satisfies the following additional condition, then $A$ is called a {\em po-semiring}:

(7) The partially ordered set $(A,\le, 0,1)$ is bounded, i.e., $1$ is the largest element and $0$ is the least element of $A$.
\end{mdef}

Note  that  condition  (7)  is  a very  strong assumption. Under the assumption,  a  po-semiring  $A$ is in fact  a  dioid, where a semiring is called a {\it dioid} if  its addition is idempotent ($ a+a=a,\,\forall a\in A$).
Furthermore, the above defined partial order $\le$ for a po-semiring $A$ is identical with the new partial order $\le_1 $ defined by the following
$$\text{$a\le_1 b$ if and only if $a+b=b$}.$$In other words, $(A,+,0,1)$ is a bounded join-semilattice.
Clearly, any bounded, distributive lattice is a po-semiring under join and meet, where $a\le b$ if and only if $a\wedge b=a$.

We remark that the class of po-semirings is much smaller than the class of semigroups. For example, the five-element lattice $M_5$ depicted in Figure 1 is  not a distributive lattice, so it is not a po-semiring. But $(M_5, \wedge)$ is clearly  a semigroup.

For a po-semiring $A$, denote by $Z(A)$ the set of all multiplicative zero-divisors. The zero-divisor graph of the multiplicative semigroup $(A,\cdot, 1)$, denoted by $\Ga(A)$, is called the zero-divisor graph of the po-semiring $A$. Clearly, all known results on zero-divisor graphs of semigroups hold for $\Ga(A)$. A nonzero element $x\in A$ is called {\em minimal}, if $0<y\leq x$ implies $x=y$ for any $y\in A$.
Note that  each minimal element of $A$ is a zero divisor, if $|Z(A)| \geq2$. Refer to \cite{Wll} for more details on po-semirings.

The prototype of a po-semiring is the po-semiring $\mathbb I(R)$ of a commutative ring $R$. The multiplication is the ideal multiplication, the addition is the addition of ideals, the partial order is the usual inclusion. Therefore, the annihilating-ideal graph of $R$ is the zero-divisor graph of the po-semiring $\mathbb I(R)$, i.e., $\mathbb{AG}(R)=\Ga(\mathbb I(R))$.

All throughout, let  $G$ be a finite or an infinite simple graph. The vertex set of $G$ is denoted by $V(G)$. The {\em core} of $G$, denoted by $C(G)$, is the largest induced subgraph of $G$ in which every edge is an edge of a cycle in $G$.
For a vertex $x$ of $G$, let $N(x)$ be the set of vertices adjacent to $x$, and call it the {\em neighborhood} of $x$. A vertex is called an {\em end vertex} if its degree is $1$. All end vertices which are adjacent to a same vertex of $G$ together with the edges is called a {\em horn}. We adopt more graph theoretic notations from \cite{BLgraph}.

We recall some notation used in \cite{LWBI}. Let $X$ and $Y$ be disjoint nonempty subsets of the vertex set of a graph. We use the notation  $X-Y$  to represent the {\em complete bipartite graph} with parts $X$ and $Y$. In particular, if $u\not\in X$, then $u-X$ represents a star graph, that is, $u$ is adjacent to every vertex in $X$ and no two distinct vertices in $X$ are adjacent. The  graph depicted in Figure 2  is called a {\em complete bipartite graph with a horn}, where the induced subgraphs on nonempty sets $X, Y, U$ are discrete. We denote the graph by $X-U-v-Y$. In particular, $G$ is called a {\em two-star graph} if $|U| =1$.
\setlength{\unitlength}{0.12cm}
\begin{figure}[h]
\centering
\input{5.LP}\nolinebreak\input{1.LP}\nolinebreak\input{3.LP}
\end{figure}

Let $\Lambda$ be an index set, and let $m,n$ be two finite or infinite cardinal numbers such that  $n=|\Lambda|\geq1$, $0\leq m\leq n $, and let $K_n$ be a complete graph with $V(K_n)=\{a_i\mid i\in \Lambda\}$. We denote the complete graph $K_n$ together with $m$ horns $X_1,X_2,\cdots, X_m$ by $K_n(m)$, where $a_1-X_1,a_2-X_2,\cdots, a_m-X_m$. For example, the graph in Figure 3 above is $K_3(2)$. Clearly, $K_1(0)$ is an isolated vertex, $K_1(1)$, $K_2(0)$ (i.e., $K_2$) and $K_2(1)$ are star graphs, $K_2(2)$ is a two-star graph, $K_3(0)$ is a triangle, while $K_3(m)$ is a triangle with $m$ horns ($1\le m\le 3$).  In this paper, we mainly study the case of $m\leq 3$. So we always assume $a_1-X,a_2-Y, a_3-Z$ for brevity.

As in \cite{LaGrange}, a graph will be called {\em realizable } (for po-semirings) if it is isomorphic to $\Ga(A)$ for some po-semiring $A$.
In this paper, we investigate the realization problem of graphs as zero-divisor graphs of po-semirings. In Section 2, it is shown that a bipartite graph $G$ is realizable if and only if
$G$ is either a complete bipartite graph or a complete bipartite graph with a horn. Realizable complete graphs with horns are then completely determined in Section 3. In particular, it is shown that $K_n(m)$ is the zero-divisor graph of some po-semiring if and only if either $0\leq m\leq min\{2,n\}$ or $m=n=3$.  The final section is devoted to the study of annihilating-ideal graph $\mathbb{AG}(R)$ of a commutative ring $R$. It is proved that $\mathbb{AG}(R)$ is a complete bipartite graph with a horn if and only if $R\cong D\times S$, where $D$ is an integral domain and $S$ is a ring with a unique non-trivial ideal. $\mathbb{AG}(R)\cong K_3(3)$ if and only if $R\cong F_1\times F_2\times F_3$, where $F_1,F_2,F_3$ are fields. We also show that there exists no ring $R$ such that $\mathbb{AG}(R)\cong K_n(2)$ for any $n\geq 3$.

Throughout the paper, set
$$
X=\{x_i\mid i\in \Gamma\},\ Y=\{y_k\mid k\in \Theta\},\ Z=\{z_p\mid p\in \Omega\},\ U=\{u_s\mid s\in \Phi\}.
$$

\section{ Bipartite graphs which are realizable for\\ po-semirings}

In this section, we give a complete classification
of all bipartite graphs which can be realized as po-semiring graphs.

\begin{lem}\label{completebipartite}
Any complete bipartite graph $G$ is realizable for po-semirings.
\end{lem}

\begin{prof}
Assume that $G$ has $X$ and $Y$ as its two parts. Set $A=\{0,1,w\}\cup X\cup Y$ and define a partial order $\le $ by
 $$
 0<x_1<x_2<\cdots<w<1,\ 0<y_1<y_2<\cdots<w<1.
 $$
 Define a commutative addition by
\begin{equation}\label{addition}
b+c=\left\{
\begin{array}{ll}
max\{b,\,c\}       & \quad\text{if $b$ and $c$ are comparable in $A$}, \\
w & \quad \text{otherwise.}
\end{array}
\right.
\end{equation}
Define a commutative multiplication on $A$ by
$$
0a=0,\ 1a=a\ (\forall a\in A),\ w^2=w,
$$
$$
x_ix_j=x_1,\ x_iy_k=0,\ x_iw=x_1,\ y_ky_l=y_1,\ y_kw=y_1\ (i,j\in \Ga, k,l\in \Theta).
$$
Clearly, $A$ is a po-semiring such that $\Ga(A)\cong G$.
\end{prof}

\begin{lem}\label{completebipartitehorn}
Let $G$ be a complete bipartite graph with a horn and set  $G: X-U-v-Y$. Then there exists a  po-semiring $A$ such that $\Ga(A)\cong G$.
\end{lem}

\begin{prof}
Let $A=\{0,1,v,w\}\cup X\cup Y\cup U$.
We define a partial order $\le $ on $A$ by
$$
0<v<x_1<x_2<\cdots<w<1,
$$
$$
0<v<y_1<y_2<\cdots<w<1,
$$
$$
0<u_1<u_2<\cdots<y_1<y_2<\cdots<w<1,
$$
Define a commutative addition by $u_i+v=y_1$ for any $u_i\in U$, and other additions are defined by (\ref{addition}).
Clearly, $(A,+)$ is a commutative semigroup.

Now we define a commutative multiplication on $A$ such that $(A, +,\cdot, 0,1)$ is a po-semiring.
For any $a\in A$, let $0a=0,\ 1a=a$. Furthermore, for any $i,j\in \Ga$, $k,l\in \Theta$ and $s,t\in \Phi$, let
$$
x_ix_j=x_iw=x_1,\  x_iy_k=x_iv=v,\  x_iu_s=0,
$$
$$
y_ky_l=y_ku_s=u_1,\  y_kv=0,\  y_kw=y_1,
$$
$$
u_su_t=u_sw=u_1,\ u_sv=v^2=0,\ vw=v,\ w^2=w.
$$
\nin A direct checking shows that the associativity holds for $(A,\cdot)$. Since the multiplication is commutative, we only need to verify that the left
distributivity holds for $(A, +,\cdot, 0,1)$. In fact, take any $a,b,c\in A$ and it is obvious,  by the addition and multiplication defined above, that $b<c$ implies $b+c=c$ and $ab\leq ac$, so $a(b+c)=ac=ab+ac$. Thus, we only need to assume that $b$ and $c$ are incomparable and show that $a(b+c)= ab+ac$.

If $a=u_s$, then
$$u_s(u_t+v)=u_sy_1=u_1=u_su_t+u_sv,$$
$$u_s(u_t+x_i)=u_sw=u_1=u_su_t+u_sx_i,$$
$$u_s(x_i+y_k)=u_sw=u_1=u_sx_i+u_sy_k.$$
Similarly, we can show that $a(b+c)=ab+ac$ always holds when $a\in\{v,x_i,y_k,w\}$ and $b,c\in A.$

Hence $(A, +,\cdot, 0,1)$ is a posemiring. Clearly, $\Ga(A)\cong G$ and the result follows.
\end{prof}

\vs{2mm}In order to check that $(A,+,\cdot, 0,1)$ defined above is a po-semiring, what is the most difficult is to check the two associative laws and one distributive law hold. In fact, each examination concerns three elements at most. By the definition of addition and multiplication in Lemma \ref{completebipartitehorn}, it is enough to take $|X| =|Y| =|U| =3$ so that $|A|=13$,  and hence we can do the examination by taking the advantage of a computer.

In general, a po-semiring corresponding to a given graph need not be unique up to isomorphism. For example, for the given complete bipartite graph with a horn $X-U-v-Y$, we can define another po-semiring $A'$ such that $A'$ is not isomorphic to the one defined in Lemma \ref{completebipartitehorn}.
\begin{ex}
Set $X+Y=\{x+y\mid x\in X,\ y\in Y\}$, and set $A'=\{0,1,v\}\cup X\cup Y \cup U\cup(X+Y)$. Let $(A',+)$ be an upper semilattice  with the Hasse diagram in Figure 4.
\setlength{\unitlength}{0.12cm}
\begin{figure}[h]\label{potwostar}
\centering
\input{2.LP}
\end{figure}
In particular, $x_1$, $y_1$ and $u_1$  are the unique minimal elements in $X$, $Y$, and $U$, respectively. Now we define a commutative multiplication on $A$ such that $(A', +,\cdot, 0,1)$ is a po-semiring. For any $i,j\in \Ga$, $k,l\in \Theta$,  $s,t\in\Phi$, we define a multiplication on $A'$ by

$$
0a=0,1a=a\ (\forall a\in A'),
$$
$$
x_ix_j=min\{x_i,x_j\}, \ x_iy_k=x_iv=v,\ x_iu_s=0,\ x_i(x_j+y_k)=x_ix_j,
$$
$$
y_ky_l=y_ku_s=u_1,\  y_kv=0,\  y_k(x_i+y_l)=y_1,\ u_su_t=u_s(x_i+y_k)=u_1,\ u_sv=0,
$$
$$
v^2=0,\ v(x_i+y_k)=v,\ (x_i+y_k)(x_j+y_l)=x_ix_j+y_1.
$$
Then $(A',+,\cdot,0,1)$ is a po-semiring with $\Ga(A')$ is isomorphic to the given graph. Clearly, $A'$ is not isomorphic to the po-semiring $A$ constructed in Lemma \ref{completebipartitehorn}, if $|X| \geq2$.
\end{ex}

\vs{3mm} The following is the main result of this section.

\begin{thm}\label{triangle-free}
Let $G$ be the zero-divisor graph of a po-semiring with $|V(G)| \geq 2$. Then the following conditions are equivalent.

{\rm(i)} $G$ is triangle-free.

{\rm(ii)} $G$ is a bipartite graph.

{\rm(iii)} $G$ is either a complete bipartite graph or a complete bipartite graph with a horn.
\end{thm}

\begin{prof}
Note that any star graph is a complete bipartite graph and any two-star graph  is a complete bipartite graph with a horn, so if $G=\Ga(S)$ for some zero-divisor semigroup $S$, then the equivalence of (i) and (ii) follows from  \cite[Theorem 2.1]{LWBI} while the equivalence of (ii) and (iii) follows from \cite[Theorem 2.10]{LWBI}. Therefore, we only need to show that both complete bipartite graphs and complete bipartite graphs with a horn can be realized for po-semiring, which follows from Lemmas \ref{completebipartite} and \ref{completebipartitehorn}.
\end{prof}

\begin{lem}\label{isolatedvertex}
Any isolated vertex $G$ is realizable for po-semirings.
\end{lem}
\begin{prof}
Let $G=\{a\}$. Set $A=\{0,1,a\}$, and define a partial order $<$ by $0<a<1$. Define a commutative addition by $x+y=max\{x,y\}$ for any $x,y\in A$. Define a commutative multiplication by
$$
0x=0,\, 1x=x \,(\forall\ x\in A),\, a^2=0.
$$
Then it is routine to check that $A$ is a po-semiring such that $\Ga(A)\cong G.$
\end{prof}

\vs{3mm} Recall that a tree is a simple connected graph $G$ without a cycle, i.e., the core of $G$ is empty. By Theorem \ref{triangle-free} and Lemma \ref{isolatedvertex} , we have the following corollary.

\begin{cor}\label{tree}
Let $G$ be a tree. Then $G$ is the graph of a po-semiring if and only if $G$ is one of the following graphs: an isolated vertex, a star graph, a two-star graph.
\end{cor}

\section{Complete graphs with horns which are realizable for po-semirings}

\vs{3mm}In this section we classify all po-semiring graphs which  are complete graphs or complete graphs with horns.
Note that the following facts were already obtained in section 2:  $K_n(m)$ is realizable for po-semirings, for any  $m,n$ with $1\leq n\leq 2$, $0\leq m\leq n$. Thus we only need to consider the case of $n\geq 3$ in this section.

\begin{lem}\label{comhorns}
For any $n\geq 3$ and $m\in\{0,1,2\}$, there exists a po-semiring $A$ such that $K_n(m)\cong\Ga(A)$.
\end{lem}
\begin{prof}
It suffices to show that each graph of the given type has a corresponding po-semiring.

{\bf Case 1.} Suppose that $m=0$. Then $G\cong K_n$ is a complete graph where $V(K_n)=\{a_i\mid i\in \Lambda\}$. Set $A=\{0,1\}\cup\{a_i\mid i\in \Lambda\}$, and define a partial order $\le$  on $A$ by $0<a_1<a_2<\cdots<1$. Define an addition on $A$ by $x+y=max\{x,y\}$ for any $x,y\in A$. Define a commutative multiplication on $A$ by
$$
0x=0,\ 1x=x\ (\forall x\in A)\ {\rm and}\ a_ia_j=0\ (\forall i,j\in \Lambda).
$$
Then it is routine to check that $A$ is a po-semiring with $\Ga(A)\cong K_n$.

{\bf Case 2.} Suppose that $m=1$. Set $A=\{0,1\}\cup\{a_i\mid i\in \Lambda\}\cup X$ and define a partial order $\le$ on $A$ by $0<a_1<a_2< \cdots <x_1<x_2<\cdots<1$. Define an addition by $x+y=max\{x,y\}$ for any $x,y\in A$. Define a commutative multiplication by
$$
0a=0,\ 1a=a\ (\forall a\in A),
$$
$$
a_ia_j=a_1x_k=0,\ a_rx_k=a_r,\ x_kx_l=x_1\ (\forall i,j,r\in \La, r\neq1, \forall k,l\in \Ga).
$$
Then $(A,+,\cdot,0,1)$ is a po-semiring such that $\Ga(A)\cong K_n(1)$, where a finite or an infinite horn $X$ is adjacent to the vertex $a_1$.

{\bf Case 3.} Suppose that $m=2$.  Denote  $X+Y=\{x+y\mid x\in X,\ y\in Y\}$, and set $A=\{0,1\}\cup\{a_i\mid i\in \La\}\cup X\cup Y \cup(X+Y)$. Let $(A,+)$ be a upper semilattice  with the Hasse diagram in Figure 5.
\setlength{\unitlength}{0.12cm}
\begin{figure}[h]
\centering
\input{4.LP}
\end{figure}

In particular, $a_n$ is the unique maximal element of $\{a_i\mid i\in \La\}$. Now we define a commutative multiplication on $A$ such that $(A, +,\cdot, 0,1)$ is a po-semiring. For any $i,j\in \La$, $k,l\in \Gamma$,  $p,q\in\Theta$, we define a multiplication on $A$ by
$$0a=0,1a=a\ (\forall a\in A),$$
$$a_ia_j=0\ (j\neq n),\ a_n^2=a_n,$$
$$a_1x_k=0,\ a_1y_p=a_1(x_k+y_p)=a_1,$$
$$a_2x_k=a_2,\ a_2y_p=0,\ a_2(x_k+y_p)=a_2,$$
$$a_nx_k=a_ny_p=a_n(x_k+y_p)=a_n,$$
$$a_ix_k=a_2,\ a_iy_p=a_1,\ a_i(x_k+y_p)=a_3,\ (i\neq 1,2,n),$$
$$x_kx_l=x_1,\ x_ky_p=a_n,\ x_k(x_l+y_p)=x_1,$$
$$y_py_q=y_p(x_k+y_q)=y_1,\ (x_k+y_p)(x_l+y_q)=x_1+y_1.$$
Then $(A,+,\cdot,0,1)$ is a po-semiring with $\Ga(A)\cong K_n(2)$.
\end{prof}

\begin{lem}\label{k33}
There exists a posemiring $A$ such that $\Ga(A)\cong K_3(3)$.
\end{lem}
\begin{prof}
Set $A=\{0,1,a_1,a_2,a_3,w\}\cup X\cup Y\cup Z$. Define a partial order on $A$ by
$$
0<a_1,a_2<z_1<z_2<\cdots<w<1,
$$
$$
0<a_2,a_3<x_1<x_2<\cdots<w<1,
$$
$$
0<a_1,a_3<y_1<y_2<\cdots<w<1.
$$
Then $A$ is a partially-ordered set with a unique maximal element $w\not=1$, where $a_1$, $a_2$ and $a_3$ are the only nonzero minimal elements. Define a commutative addition by
$$
0+a=a+a=a,\ 1+a=1,\ (\forall a\in A),
$$
$$
\ a_1+a_2=z_1,\ a_1+a_3=y_1,\ a_2+a_3=x_1,
$$
and for any $b,c\in A$ such that $\{b,c\}\nsubseteq \{a_1,a_2,a_3\}$, we define
\begin{equation}
b+c=\left\{
\begin{array}{ll}
max\{b,\,c\}       & \quad\text{if $b$ and $c$ are comparable in $A$}, \\
w & \quad \text{otherwise.}
\end{array}
\right.
\end{equation}
For any $k\in \Gamma$, $p\in \Theta$ and $s\in\Omega$, define a commutative multiplication by
$$
0a=0,\ 1a=a\ (\forall a\in A),
$$
$$
a_1^2=a_1, a_1a_2=a_1a_3=a_1x_k=0, a_1y_p=a_1z_s=a_1w=a_1,
$$
$$
a_2^2=a_2, a_2a_3=a_2y_p=0, a_2x_k=a_2z_s=a_2w=a_2,
$$
$$
a_3^2=a_3, a_3z_s=0, a_3x_k=a_3y_p=a_3w=a_3,
$$
$$
X^2=\{x_1\}, XY=\{a_3\}, XZ=\{a_2\},wX=\{x_1\},
$$
$$
Y^2=\{y_1\},YZ=\{a_1\}, wY=\{y_1\}, Z^2=\{z_1\}, wZ=\{z_1\}, w^2=w.
$$
It is easy to check that $(A,+,\cdot,0,1)$ is a po-semiring such that $\Ga(A)\cong K_3(3)$.
\end{prof}

\vs{3mm}The following is an improvement of \cite[Theorem 2.2]{WL2}.

\begin{lem}\label{nleqmleq4}
For any $n\geq4$ and any $m$ with $n\geq m\geq3$, $K_n(m)$ has no corresponding semigroups.
\end{lem}
\begin{prof}
Assume to the contrary that $S=\{0\}\cup \{a_i\mid i\in \La\}\cup X_1\cup X_2\cup\cdots \cup X_m$ is a zero-divisor semigroup such that $\Ga(S)\cong K_n(m)$, where $n\geq4$ and $n\geq m\geq3$.

First we claim that for any  $x_j\in X_j$, we must have $a_ix_j=a_i$,  where $1\leq i\neq j\leq m$. In fact, if $a_ix_j=x_k$ for some $x_k\in X_k$, then $a_jx_k=a_i(a_jx_j)=0$
so that $j=k$ whence $x_ix_j=x_i(a_i x_j)=0$, a contradiction. Thus, $a_ix_j\in \{a_i\mid i\in \La\}$. Let $a_ix_j=a_r$, then for any $x_i'\in X_i$, $a_rx_i'=(a_ix_i')x_j=0$ and hence $r=i$, that is, $a_ix_j=a_i.$

Take $x_2\in X_2$, $x_3\in X_3$. Then $a_1(x_2x_3)=(a_1x_2)x_3=a_1x_3=a_1$ so that $x_2x_3\not\in\{a_i\mid i\in \La, i\neq 1\}$.  On the other hand, $a_2(x_2x_3)=0=a_3(x_2x_3)$ which means that $x_2x_3\in\{a_i\mid i\in \La\}$ and hence $x_2x_3=a_1$. So $(a_4x_2)x_3=a_4a_1=0$, which implies that $a_4x_2=a_3$. Substituting $a_3,x_3,a_1,x_1$ by $a_1,x_1,a_3,x_3$ respectively, we can obtain that $a_4x_2=a_1$, a contradiction.
\end{prof}

\vs{2mm}Note that a po-semiring is certainly a multiplicative commutative semigroup. So we have
\begin{cor}
For any $n\geq4$ and any $m$ with $n\geq m\geq3$, $K_n(m)$ has no corresponding po-semirings.
\end{cor}

In view of Lemmas \ref{isolatedvertex}, \ref{comhorns}, \ref{k33} and \ref{nleqmleq4}, we obtain the main result of this section.
\begin{thm}\label{completehorns}
Let $G=K_n(m)$ be a complete graph $K_n$ together with $m$ horns, where $n$ and $m$ could be any cardinal numbers such that $0\leq m\leq n$ and $n\geq1$. Then $G$ is a po-semiring graph if and only if either $0\le m\le min\{2,n\}$ or $n=m=3$.
\end{thm}

\section{ Some results on annihilating-ideal graphs of commutative rings }

The annihilating-ideal graph $\mathbb{AG}(R) $ of a commutative ring $R$ was introduced in \cite{BR} and further studied in \cite{BRII, AA2,AA,AANJS}. In this section, we add more results on the graph.

\begin{lem}\label{minimal}
Let $A$ be a po-semiring and let $x-u-y$ be a path in $\Ga(A)$. If it is contained in neither triangle nor quadrilateral, then $u$ is a minimal
element of $A$.
\end{lem}

\begin{prof}
 Assume to the contrary that $u$ is not a minimal element of $A$. Assume further $0<v<u$ for some $v\in A$. If $v=x,$ then $xy=vy\leq uy=0$ whence there is a triangle $x-u-y-x$, a contradiction. Similarly, $v\neq y$.
If $v\not=x,y$, then there is a square $x-u-y-v-x$, another contradiction.
\end{prof}

\begin{lem}\label{maximalideal}
Let $R$ be a ring, $I$ a minimal ideal of $R$. Then $ann(I)$ is a maximal ideal of $R$.
\end{lem}
\begin{prof}
Since $I$ is minimal, $I$ is a principal ideal. Suppose that $I=Rx$ for some $x\in I$. Then $ann(x)=ann(I)$ whence $Rx\cong R/ann(I)$ is a simple $R$-module, so $ann(I)$ is a maximal ideal.
\end{prof}

\vs{3mm}In view of \cite[Theorem 2.3]{BRII} and \cite[Corollary 23]{AANJS}, we known that $\mathbb{AG}(R)$ is a complete bipartite graph if and only if either $\mathbb{AG}(R)$ is a star graph or $R$ is a reduced ring with $|Min(R)| =2$. For a complete bipartite graph with a horn, we have

\begin{mthm}\label{annihilatingidealcompletebiaprtitewithhorn}
The graph $\mathbb {AG}(R)$ is a complete bipartite graph with a horn if and only if $R\cong D\times S$, where $D$ is an integral domain and $S$ is a ring with a unique non-trivial ideal. In this case, $\mathbb {AG}(R)\cong K_1+D_r+K_1+D_r$ where $r$ is either $1$ or an infinite cardinal number.
\end{mthm}
\begin{prof}
$($$\Rightarrow$$)$ Let $\mathbb{AG}(R)$ be  a complete bipartite graph with a horn£º $X-U-\mathfrak{c}-Y$. Set $\mathfrak{p}=\sum_{\mathfrak p_i\in X}\mathfrak{p}_i$, $\mathfrak{q}=\sum_{\mathfrak q_j\in Y}\mathfrak{q}_j$ and $\mathfrak{u}=\sum_{\mathfrak u_k\in U}\mathfrak{u}_k$. Clearly, $\mathfrak {cq}=0$. But for any $\mathfrak{p}_i\in X,\mathfrak{u}_k\in U$, $\mathfrak{qp}_i\neq0,\mathfrak{qu}_k\neq0$, which means that $\mathfrak q\in Y$. Similarly, we have $\mathfrak p\in X$ and $\mathfrak u\in U$.

If $|U| =1$, then $\mathbb {AG}(R)$ is a two-star graph. By \cite [Theorem 2]{AA}, $\mathbb {AG}(R)\cong P_4$ and  $R\cong F\times S$ where $F$ is a field and $S$ is a ring with a unique non-trivial ideal. So we assume that $|U| \geq2$ in next discussion.

 Now, we show the following claims:

{\bf Claim 1.} $\mathfrak{q}$ is a maximal ideal. By Lemma \ref{minimal}, $\mathfrak{c}$ is a minimal ideal. Note that $\mathfrak {cq}=0$, so we only need to show that $ann(\mathfrak c)=\mathfrak{q}$ by Lemma \ref{maximalideal}. In fact, first we have $\mathfrak{p}\mathfrak{q}\neq 0$. On the other hand, for any $\mathfrak{u}_j\in U$, $\mathfrak u_j(\mathfrak{p}\mathfrak{q})=0=\mathfrak{c}(\mathfrak{p}\mathfrak{q})$ which together with $|U| \geq2$ yield that $\mathfrak{p}\mathfrak{q}=\mathfrak{c}$. Thus, $\mathfrak{c}^2=(\mathfrak{c}\mathfrak{q})\mathfrak{p}=0$ whence $\mathfrak{c}(\mathfrak{c}+\mathfrak{q})=0$, which together with $\mathfrak{u}(\mathfrak{c}+\mathfrak{q})=\mathfrak{u}\mathfrak{q}\neq0$ and $\mathfrak{p}(\mathfrak{c}+\mathfrak{q})\neq0$ imply that $\mathfrak{c}+\mathfrak{q}\in Y$ and hence $\mathfrak{c}\seq \mathfrak{q}$. Since $\mathfrak{c}(\mathfrak{u}+\mathfrak{q})=0$, but $\mathfrak{p}(\mathfrak{u}+\mathfrak{q})=\mathfrak{pq}=\mathfrak c\neq0$ so that $\mathfrak{u}+\mathfrak{q}\in Y $, which means that $\mathfrak{u}\seq \mathfrak{q}$. Therefore, we obtain that $ann(\mathfrak{c})=\mathfrak c\cup\mathfrak u\cup\mathfrak{q}=\mathfrak{q}$.

{\bf Claim 2.} $\mathfrak{p}$ is a prime ideal and $\mathfrak{p}+\mathfrak{q}^2=R$. Take any  $x\in \mathfrak{u}$. Then $\mathfrak p(Rx)\seq \mathfrak {pu}=0=\mathfrak c(Rx)$, so $Rx\in U$. Clearly, $x^2\neq 0$. Otherwise, $(Rx)(Rx+\mathfrak c)=0=\mathfrak c(Rx+\mathfrak c)$. Since $\mathfrak c$ is minimal, it follows that $Rx+\mathfrak c\neq \mathfrak c$, so $Rx+\mathfrak c\in U$ and hence $Rx+\mathfrak c=Rx$. Thus $\mathfrak c\seq Rx$ and hence $\mathfrak {cp}\seq \mathfrak p(Rx)=0$, a contradiction. So $x^2\neq 0$ whence $ann(x)=\mathfrak c\cup\mathfrak p$. Since $\mathfrak c$ is minimal, it follows that $\mathfrak {cp}=\mathfrak c$ whence $\mathfrak c\seq \mathfrak p$ and hence $ann(x)=\mathfrak p$. If $\mathfrak p$ is not a maximal ideal contained in $Z(R)$, then there exists $y\in Z(R)\backslash \mathfrak p$ such that $\mathfrak p+Ry\seq Z(R)$. By the proof of Claim 1, we obtain that $\mathfrak u \seq \mathfrak q$ and $\mathfrak c\seq \mathfrak q$, so we have $Z(R)=\mathfrak p\cup \mathfrak u \cup \mathfrak c\cup \mathfrak q=\mathfrak p\cup \mathfrak q$ whence $y\in \mathfrak q\backslash \mathfrak p$. Take any $z\in \mathfrak p$ such that $\mathfrak cz\neq 0$. Then $Rz\in X$ and hence $\mathfrak c(y+z)=\mathfrak cz\neq 0$ so that $y+z\not\in \mathfrak q$, at the same time, $\mathfrak u(y+z)=\mathfrak uy\neq0$ which implies that $y+z\not\in \mathfrak p$, that is, $y+z\not\in Z(R)$, a contradiction. This proves that $\mathfrak p$ is maximal among all ideals of $R$ that are annihilators of elements, so $\mathfrak p$ is a prime ideal. Clearly, $\mathfrak p$ and $\mathfrak q$ are incomparable. Otherwise, we have $\mathfrak p\seq \mathfrak q$, so $\mathfrak c=\mathfrak {pc}\seq \mathfrak{qc}=0$, a contradiction. Therefore, $\mathfrak{p}+\mathfrak{q}^2=R$, as required.

{\bf Claim 3.}  $\mathfrak p\cap \mathfrak q^2=0$. Note that for any $\mathfrak u_k\in U$, we have $\mathfrak u_k(\mathfrak p\cap \mathfrak q^2)\seq \mathfrak u_k\mathfrak p=0$ and $\mathfrak c(\mathfrak p\cap \mathfrak q^2)\seq\mathfrak c\mathfrak q =0$ which together with $|U| \geq 2$ yield  that $\mathfrak p\cap \mathfrak q^2=\mathfrak c$ or $\mathfrak p\cap \mathfrak q^2=0$. If $\mathfrak p\cap \mathfrak q^2=\mathfrak c$, then $\mathfrak c\seq \mathfrak q^2$, hence $\mathfrak c=\mathfrak {pc}\seq \mathfrak {pq}^2=\mathfrak {cq}=0$, a contradiction. So, $\mathfrak p\cap \mathfrak q^2=0$.

By Chinese Remainder Theorem, we have $R\cong R/\mathfrak p\times R/\mathfrak q^2$. Clearly, $R/\mathfrak p$ is an integral domain, $R/\mathfrak q^2$ is an Artinian local ring. Thus all non-trivial ideals of $R/\mathfrak q^2$ are annihilating-ideals. If there exist two non-trivial ideals, say $\mathfrak{m,n}$, in
$R/\mathfrak q^2$, then $\mathfrak{m,n}\seq \mathfrak q/\mathfrak q^2$ and hence $\mathfrak{mn}=0$. So we have a triangle $(R/\mathfrak p,0)-(0,\mathfrak m)-(0,\mathfrak n)-(R/\mathfrak p,0)$ in $\mathbb {AG}(R)$, which contradicts the fact that $\mathbb {AG}(R)$ is a bipartite graph. This complete the proof of necessity.

$($$\Leftarrow$$)$ Let $R=D\times S$ where $D$ is an integral domain, $S$ has a unique non-trivial ideal $\mathfrak m$. If $D$ is not a field, then $D$ has infinitely many ideals, and $\mathbb {AG}(R)$ is the complete bipartite graph with a horn $X-U-\mathfrak c-Y$, where
$$
\mathfrak c=(0,\mathfrak m),\ X=\{(0,S)\},\\
$$
$$
Y=\{(\mathfrak b,\mathfrak m)\mid \mathfrak b\in \mathbb I(D)^*\},\
U=\{(\mathfrak b,0)\mid \mathfrak b\in \mathbb I(D)^*\}.
$$
So $\mathbb{AG}(R)\cong K_1+D_r+K_1+D_r$, where $r$ is the  cardinality of $\mathbb I(D)^*$. In this case, $r=\infty$.
If $D$ is a field, then $r=1$ and $\mathbb {AG}(R)\cong P_4$.  This complete the proof of the theorem.
\end{prof}

\begin{rem}\label{onlyoneidealring}
Let $R$ be a ring with only one non-trivial ideal. Then, by \cite{JR}, either $R\cong K[x]/(x^2)$ where $K$ is a field or $R\cong V/p^2V$, where $V$ is a discrete valuation ring of {\it characteristic zero} and residue field of {\it characteristic} $p$, for some prime number $p$. In the latter case, the characteristic of $R$ is $p^2$.
\end{rem}

By \cite[Theorem 3]{AA2}, the graph $\mathbb{AG}(R)$ is a complete graph if and only if either $\mathbb{AG}(R)\cong K_2$ or $Z(R)$ is an ideal of $R$ with $Z(R)^2=0$. Moreover, in the first case, either $R\cong F_1\times F_2$, where $F_1,F_2$ are fields, or $(R, \mathfrak m)$ is a local ring with exactly two non-trivial ideals $\mathfrak m$ and $\mathfrak m^2$.
For complete graphs with horns, $K_1(0)$ is an isolated vertex, $K_1(1), K_2(0)$ and $K_2(1)$ are star graphs, while $K_2(2)$ is a two-star graph. By \cite[Theorem 2]{AA2} or Theorem \ref{annihilatingidealcompletebiaprtitewithhorn} above, $\mathbb{AG}(R)$ is a two-star graph if and only if $\mathbb{AG}(R)\cong P_4$ , if and only if $R\cong F\times S$, where $F$ is a field and $S$ is a ring with a unique non-trivial ideal. For $K_n(m)$ where $n\geq 3$ and $1\leq m\leq 3$, we have the following results.

\begin{lem}
Let $\mathbb{AG}(R)\cong K_n(1)$, where $n\geq3$. Then $Z(R)$ is a maximal ideal. If further $R$ is  Artinian, then $Z(R)^5=0$.
\end{lem}
\begin{prof}
Let $\mathfrak c$ be the only center of $\mathbb{AG}(R)$, that is, the only horn is adjacent to $\mathfrak c$. Then, by Lemma \ref{minimal}, $\mathfrak c$ is a minimal ideal. Take distinct $\mathfrak {a,b}\in V(K_n)\backslash \{\mathfrak c\}$. Then $\mathfrak a(\mathfrak c+\mathfrak b)=0$ so that $\mathfrak c+\mathfrak b\in V(K_n)$. Obviously, $\mathfrak c+\mathfrak b\neq \mathfrak c$ and hence $\mathfrak c(\mathfrak c+\mathfrak b)=0$, that is, $\mathfrak c^2=0$. Put $\mathfrak{m}=Z(R)$. Clearly, $\mathfrak{m}=ann(\mathfrak c)$, which means that $\mathfrak m$ is maximal by Lemma \ref{maximalideal}. If $R$ is an Artinian ring, then  $\mathfrak{m}$ is a nilpotent ideal. If the nilpotency index is $n\geq6$, then $\{\mathfrak{m}^{n-2},\mathfrak{m}^{n-1}\}\seq N(\mathfrak{m}^2)\cap N(\mathfrak{m}^3)$, so $\mathfrak{m}^2$ and $\mathfrak{m}^3$ can not be end vertices, Hence $\mathfrak{m}^2,\mathfrak{m}^3\in V(K_n)$ so that
$\mathfrak{m}^5=0$, a contradiction. Therefore, $n\leq 5$, as required.
\end{prof}

\vs{3mm}The structure of $R$ with $\mathbb{AG}(R)\cong K_n(1)$  seems to be rather complicated and hard to determine completely. For $n=1$ or $n=2$ where the unique horn contains exactly one vertex,  see \cite{WULUoneortwo} for the detailed structure theorems on commutative rings with at most three nontrivial ideals. We will discuss the problem for $n\ge 3$ in a separate paper.
\begin{mthm}\label{annikn2}
There exists no ring $R$ such that $\mathbb{AG}(R)\cong K_n(2)$ for any $n\geq 3$.
\end{mthm}
\begin{prof}
Let $R$ be a ring such that $\mathbb{AG}(R)\cong K_n(2)$ for some $n\geq 3$.
Suppose that $K_n$ is the complete graph with $V(K_n)=\{a_i\mid i\in \Lambda\}$ where $n=|\La| \geq 3$, and the two horns are $\mathfrak a_1-X$, $\mathfrak a_2-Y$. Put $\mathfrak a=\sum_{i\in \La}\mathfrak a_i$, $\mathfrak p=\sum_{\mathfrak p_j\in X}\mathfrak p_j$ and $\mathfrak q=\sum_{\mathfrak q_k\in Y}\mathfrak q_k$. Then $\mathfrak p\in X$, $\mathfrak q\in Y$. Now, we have the following claims:

{\bf Claim 1.} $\mathfrak p $ and $\mathfrak q$ are maximal ideals.
 In view of Lemma \ref{minimal}, we get that $\mathfrak a_1,\mathfrak a_2$ are minimal, so $\mathfrak a_1+\mathfrak a_2\not\in \{\mathfrak a_1,\mathfrak a_2\}$. Since $\mathfrak a_3(\mathfrak a_1+\mathfrak a_2)=0$, it follows that $\mathfrak a_1+\mathfrak a_2\in \{\mathfrak a_i\mid i\in \La, i\neq 1,2\}$ and hence $\mathfrak a_1(\mathfrak a_1+\mathfrak a_2)=0$, so $\mathfrak a_1^2=0$, which yields that $\mathfrak a_1\mathfrak a=0$. Similarly, $\mathfrak a_2^2=0, \mathfrak a_2\mathfrak a=0$, which means that $\mathfrak a\in \{\mathfrak a_i\mid i\in \La\}$. Note that $\mathfrak a_1(\mathfrak a+\mathfrak p)=0$, it follows that $\mathfrak a+\mathfrak p\in \mathbb {A}(R)$. On the other hand, $\mathfrak a_2(\mathfrak a+\mathfrak p)=\mathfrak a_2\mathfrak p\neq 0$, so $\mathfrak a+\mathfrak p\in X$ and hence $\mathfrak a+\mathfrak p\seq \mathfrak p$ whence $\mathfrak a\seq \mathfrak p$. Thus $ann(\mathfrak a_1)=\mathfrak a\cup\mathfrak p=\mathfrak p$. Consequently, $\mathfrak p$ is a maximal ideal of $R$ by Lemma \ref{maximalideal}. Similarly, $\mathfrak q$ is also maximal.

{\bf Claim 2.} $\mathfrak a^3=0$.  By the proof of Claim 1, $\mathfrak a\seq\mathfrak p\cap\mathfrak q$. On the other hand, $\mathfrak a_1(\mathfrak p\cap\mathfrak q)=\mathfrak a_2(\mathfrak p\cap\mathfrak q)=0$ implies that $\mathfrak p\cap\mathfrak q\in V(K_n)$ and hence $\mathfrak p\cap\mathfrak q\seq \mathfrak a$ so that $\mathfrak p\cap\mathfrak q= \mathfrak a$. Now, we show that $\mathfrak a^2\neq\mathfrak a$.  Assume to the contrary that $\mathfrak a^2=\mathfrak a$. Since $\mathfrak a_1+\mathfrak a_2\seq \mathfrak a$, it follows that $\mathfrak a$ is not minimal. So we can take $0\neq x\in\mathfrak a$ such that $\mathfrak a\neq Rx$. Then there exist $u,v\in \mathfrak a$ such that $x=uv$. This implies that $(Ru)(Rv)=Rx\neq 0$. If $\mathfrak a=Ru=Rv$, then $Rx=\mathfrak a^2=\mathfrak a$, a contradiction. If at least one of $Ru,Rv$ is not equal to $\mathfrak a$, without loss of generality, suppose that $Ru\neq \mathfrak a$, then $Ru\in V(K_n)$ whence $(Ru)(Rv)\seq (Ru)\mathfrak a=0$, that is, $Rx=0$, another contradiction. Thus, $\mathfrak a^2\neq \mathfrak a$. Clearly, $\mathfrak a^2\in V(K_n)$, hence $\mathfrak a^3=\mathfrak a\mathfrak a^2=0$.

By the proof of Claim 2, $\mathfrak p\cap\mathfrak q=\mathfrak a$, so $\mathfrak p\mathfrak q\seq \mathfrak a$ and hence $\mathfrak p^3\mathfrak q^3=0$. By Claim 1, $\mathfrak p,\mathfrak q$ are maximal, hence $R$ is an Artinian ring with exactly two maximal ideals. By \cite [Theorem 8.7]{AM}, there exist two Artinian local rings, say $(R_1,\mathfrak m)$ and $(R_2,\mathfrak n)$, such that $R=R_1\times R_2$. Hence we can assume that $\mathfrak p=(\mathfrak m,R_2)$, $\mathfrak q=(R_1,\mathfrak n)$ and hence $\mathfrak a=(\mathfrak m,\mathfrak n)$. Note that $(R_1,0)-(0,R_2)-(\mathfrak m,0)$ is a path in $\mathbb {AG}(R)$, it follows that $(0,R_2)$ is not an end vertices, that is, $(0,R_2)\in V(K_n)$ and hence $\mathfrak a(0,R_2)=(0,\mathfrak n)=0$ so that $\mathfrak n=0$. Similarly, $\mathfrak m=0$. This implies that $R_1$ and $R_2$ are fields, and $\mathbb {AG}(R)\cong K_2$, a contradiction. The contradiction followed from the assumption that there exists a ring $R$ such that $\mathbb{AG}(R)\cong K_n(2)$ for some $n\geq 3$.
\end{prof}

\begin{mthm}\label{annik33}
$\mathbb{AG}(R)$ is $K_3(3)$ if and only if $R\cong F_1\times F_2\times F_3$, where $F_1,F_2,F_3$ are fields.
\end{mthm}

\begin{prof}
Let $\mathbb{AG}(R)\cong K_3(3)$. Suppose that $V(K_3)=\{\mathfrak {a,b,c}\}$,  the three horns are $\mathfrak a-X$, $\mathfrak b-Y$ and $\mathfrak c-Z$. By Lemma \ref{minimal}, $\mathfrak {a,b,c}$ are minimal ideals. Denote $\mathfrak m=\sum_{\mathfrak m_i\in X}\mathfrak m_i$, $\mathfrak n=\sum_{\mathfrak n_j\in Y}\mathfrak n_j$, $\mathfrak p=\sum_{\mathfrak p_k\in Z}\mathfrak p_k$. Clearly, $\mathfrak {am}=0$, but $\mathfrak {bm}\neq0$, $\mathfrak {cm}\neq0$, so $\mathfrak m\in X$. Similarly, $\mathfrak n\in Y$, $\mathfrak p\in Z$. Now, we show the following claims:

{\bf Claim 1.} $\mathfrak a^2=\mathfrak a,\mathfrak b^2=\mathfrak b,\mathfrak c^2=\mathfrak c$. Since $\mathfrak a$ is minimal, it follows that $\mathfrak a^2\in \{0,\mathfrak a\}$. If $\mathfrak a^2=0$, then $\mathfrak a(\mathfrak a+\mathfrak b)=0$, which together with $\mathfrak c(\mathfrak a+\mathfrak b)=0$ yields that $\mathfrak a+\mathfrak b\in \{\mathfrak {a,b,c}\}$, which contradicts the fact $\mathfrak {a,b,c}$ are minimal ideals. Similarly, $\mathfrak b^2=\mathfrak b,\mathfrak c^2=\mathfrak c$.

{\bf Claim 2.} $\mathfrak {m,n,p}$ are maximal ideal of $R$. Since $\mathfrak a\mathfrak m=0$, so, by Lemma \ref{maximalideal}, we only need to show that $\mathfrak m=ann(\mathfrak a)$. Note that $\mathfrak a(\mathfrak m+\mathfrak b)=0$, $\mathfrak b(\mathfrak m+\mathfrak b)\neq0$, $\mathfrak c(\mathfrak m+\mathfrak b)\neq0$, hence $\mathfrak m+\mathfrak b\in X$ whence $\mathfrak m+\mathfrak b\seq \mathfrak m$, that is, $\mathfrak b\seq\mathfrak m$. Similarly, $\mathfrak c\seq\mathfrak m$. Therefore, $ann(\mathfrak a)=\mathfrak b\cup\mathfrak c\cup\mathfrak m=\mathfrak m$ so that $\mathfrak m$ is a maximal ideal. The other two results can be obtained similarly.

{\bf Claim 3.} $\mathfrak m\cap\mathfrak n\cap\mathfrak p=0$. Take any $ x\in \mathfrak m\cap\mathfrak n\cap\mathfrak p$. If $x\neq 0$, then $Rx\in N(\mathfrak a)\cap N(\mathfrak b)\cap N(\mathfrak c)=\emptyset$, a contradiction.

By Chinese Remainder Theorem, we have $R\cong R/\mathfrak m\times R/\mathfrak n\times R/\mathfrak p$, that is,
$R\cong F_1\times F_2\times F_3$, where $F_1,F_2,F_3$ are fields. The converse is trivial.
\end{prof}

\begin{cor}
If $\mathbb{AG}(R)$ is $K_3(3)$, then $|\mathbb{AG}(R)| =6$.
\end{cor}

It is interesting to make a comparison between the results of sections 2, 3 and those in section 4. Lemma \ref{completebipartitehorn} shows that for any complete bipartite graph with a horn, there exists a corresponding  po-semiring. However, in view of Theorem \ref{annihilatingidealcompletebiaprtitewithhorn}, only a few of complete bipartite graphs with a horn can ba realized as annihilating-ideal graph. We can obtain a similar result for complete graphs with horns by Theorems 3.5, 4.6 and 4.7.
 %\ref{completehorns} and Theorem \ref{annikn2} and \ref{annik33}.
This observation indicates that  the class of po-semirings $\mathbb I(R)$ of  rings $R$ is a very small subclass of po-semirings. So, we have the following question.

\begin{que}
Classify the po-semirings $A$ such that $A\cong\mathbb I(R)$ for some commutative ring $R$.
\end{que}

\end{document}